\documentclass{amsart}
\usepackage[T2A]{fontenc}
\usepackage[russian,english]{babel}
\usepackage[a4paper,textwidth=166mm,textheight=242mm,centering]{geometry}

\def\a{\alpha}

\newtheorem{теор}{Тheorem}[section]
\newtheorem{лем}[теор]{Lemma}
\newtheorem{зам}[теор]{Remark}

\newtheorem{анн}[теор]{Annotation}
\newtheorem{опр}[теор]{Definition}
\newtheorem{ключ}[теор]{ Keywords.}

\newtheorem{thm}{Theorem}[section]

\begin{document}

\begin{center}
\textbf{{\large {\ 
Space--Time Source Recovery for a Subdiffusion Equation from Integral Observations}}}\\[0pt]
\medskip \textbf{R.R. Ashurov$^{1}$ and O.T. Mukhiddinova$^{1,2}$}\\[0pt]
\textit{ashurovr@gmail.com, oqila1992@mail.ru \\[0pt]}

\smallskip
\textit{{$^{1}$V.I. Romanovskiy Institute of Mathematics, Uzbekistan Academy of Science, University str.,9, Olmazor district, Tashkent, 100174, Uzbekistan}\\
{$^{2}$ Tashkent University of Information Technologies named after Muhammad al-Khwarizmi,
Str., 108, Amir Temur Avenue,Tashkent, 100200, Uzbekistan } }\\

\end{center}

\begin{анн}
This paper is devoted to the study of the inverse problem of determining the right-hand side of the subdiffusion equation with the Caputo derivative with respect to time. In our case, the inverse problem consists in restoring the coefficient of the right-hand side, which depends on both the time and the spatial variable, when measured in integral form. Previously, similar inverse problems were studied for hyperbolic and parabolic equations with a different overdetermination condition, and in some works the existence and uniqueness of generalized solutions was established, while in others, the uniqueness of classical solutions was established. However, similar inverse problems for fractional equations with an integral overdetermination condition have not been considered before this work. The existence and uniqueness of a weak solution to the inverse problem under consideration is established. It is noteworthy that the results obtained are new for parabolic equations as well.

\end{анн}

\begin{ключ}
Subdiffusion equation, Caputo time derivative, inverse problem, uniqueness and existence of solution, Fourier method.
\end{ключ}
\section{Introduction}

In recent decades, much attention has been paid to the study of problems for partial differential equations
with fractional derivatives, which are used to model a number of phenomena in various branches of science
and technology (see, e.g., \cite{Machado}–\cite{SU}). It should also be noted that, according to some recent research results, it can be seen that there are some physical phenomena that cannot be modeled by partial differential equations
of integer order, but are adequately described by fractional derivatives (\cite{Herrmann} - \cite{Podlubny}). On the other hand, partial differential equations with fractional derivatives, along with such important physical applications, in some cases
generalize partial differential equations of integer order (see, e.g., \cite{Kil}–\cite{ACT}).
Therefore, the study of direct and inverse problems associated with differential equations with fractional
derivatives is of great interest from both theoretical and practical points of view.

In recent years, there has been growing interest among researchers in inverse problems of determining the source function in differential equations of both integer and fractional order (see, for example, the monograph by S. I. Kabanikhin \cite{Kaban} and the review article by Yamamoto \cite{Yamamoto1}). This interest is due to their crucial role in applications in various fields, such as mechanics, seismology, medical tomography, and geophysics (see, for example, \cite{Machado}).
Mainly such problems have been studied where the source function has the form \( F(x,t) = g(t)f(x) \), with either \( g(t) \) or \( f(x) \) being unknown. To our knowledge, the more general case \( F(x,t) \) without such factorization has not been studied. In this scenario, even the choice of an appropriate overdetermination condition remains unclear. Inverse problems aimed at recovering the time-dependent component \( g(t) \) are usually solved by reducing them to integral equations (see, e.g., \cite{ASh1}–\cite{ASh2} and references therein). In contrast, the inverse problem of determining the space-dependent component \( f(x) \) is analyzed in two different cases: \( g(t) \equiv 1 \) and \( g(t) \not\equiv 1 \). When \( g(t) \equiv 1 \), such problems were studied, e.g., in \cite{KirM}–\cite{AMux1}. The case when \( g(t) \not\equiv 1 \) is more complicated, and solvability depends on the sign-definiteness of the function \( g(t) \) (see, for example, \cite{ASh3}–\cite{VanBockstal}, as well as \cite{Kaban} and the review article \cite{Yamamoto1}).

The works most closely related to our study are those by S. Z. Dzhamalov et al. \cite{Djamalov1}–\cite{Djamalov5}, which address inverse problems of determining the right-hand side of equations, where the right-hand side depends on both time and a portion of the spatial variable. These studies investigate the existence of a generalized solution to the problem using the Galerkin method.

In \cite{Fikret}, the inverse problem for a hyperbolic equation is also addressed, specifically the determination of the right-hand side, which depends on both time and a portion of the spatial variable. The authors successfully established the uniqueness of the classical solution to this inverse problem.

In recent work \cite{AshurovMuhiddinovanew} the authors investigated the inverse problem of determining the right-hand side of a subdiffusion equation with a Caputo time derivative, where the right-hand side depends on both time and certain spatial variables.   The primary objective of this paper is to investigate the inverse problem of determining the coefficient of the right-hand side, and solution of the forward problem: \( \{ u(t, x, y), h(t, x) \} \), \((x,\, y)\, \in (0,1) \times (0, \pi) \subset \mathbb{R}^2\), subject to the overdetermination condition given by:
\begin{equation}\label{l_0}
    u(t, x, l_0) = \psi(t, x), \quad t \in (0, T), \quad x \in (0, 1), \quad l_0 \in (0, \pi),
\end{equation}
where \( \psi \) is a known continuous function.
 It is proved the existence and uniqueness of the weak solution to the considered inverse problem. To solve it, the authors employed the Fourier method with respect to the variable independent of the unknown right-hand side, followed by the method of successive approximations to compute the Fourier coefficients of the solution. 

 Note that in all works \cite{Djamalov1} - \cite{{AshurovMuhiddinovanew}} the authors considered the inverse problem with the overdetermination condition (\ref{l_0}).

In this paper, we set the overedetermination condition in  a integral form and investigate the inverse problem of determining the function \( h(t, x) \) that appears in the right-hand side of a subdiffusion equation, given as \( f(t, x, y) h(t, x) + g(t, x, y) \), where \( x \in \mathbb{R}^m \) and \( y \in \mathbb{R}^n \). 

Next, we proceed to the precise formulation of the problem.

 Let $ Q=(0,T) \times G_x\times  \Omega_y$
and
$ D= G_x \times \Omega_y$, where $G:= G_x\subset \mathbb{R}^m$ and $\Omega:= \Omega_y \subset \mathbb{R}^n$ are a bounded domains, with a sufficiently smooth boundary.
Consider the following initial - boundary value problem
\begin{equation} \label{1}
	\begin{cases}
		 D_t^{\alpha}u - \Delta_x u - \Delta_y u = g(t,x,y) + f(t,x,y) \cdot h(t,x) , \quad (t,x,y) \in Q, \\
        
         u(0,x,y) = \varphi(x,y), \quad (x,y) \in D, \\
         
         u(t,0,y) =  u(t,1,y) = 0,\quad t\in [0,T], \quad y\in \Omega,\\
         
         u(t,x,0) = u(t,x,\pi) = 0,\quad t\in [0,T], \quad x\in G.
	\end{cases}
\end{equation} 
 Here $\Delta_x$ and $\Delta_y$ are the Laplace operators in variables $x$ and $y$ respectively, \( f, g \) and \( \varphi \) are given functions, \( \alpha \in (0,1) \), and \( D_t^{\alpha} \) is the fractional Caputo derivative. Recall, if \( v(t) \) is  an absolutely continuous function, then its Caputo derivative has the form (see, e.g., \cite{Kil}, p. 91)
\[
D_{t}^{\alpha } v(t)\equiv J_t^{1-\alpha} D_t v(t), \quad J_t^\a v(t) = \frac{1}{\Gamma (\alpha )} \int _{0}^{t} (t-\tau )^{\alpha-1 }v(\tau ) d\tau ,
\] 
where $D_t=d/dt$, $\Gamma(\alpha)$ is the gamma function, $J_t^\alpha$ is the Riemann-Liouville fractional integral.

If the function \(h(t, x)\) is given, then problem (\ref{1}) is called \textit{the forward problem} and its solution exists and is unique under certain conditions on the problem’s data (see, e.g., \cite{Sakamoto}, \cite{AMux2}).

Now we will assume that the coefficient \( h(t, x) \) of the source function is unknown and must be determined. The main goal of this paper is to study \textit{the inverse problem }of determining a pair of functions \( \{ u(t, x, y), h(t, x) \} \) under the overdetermination condition specified by the following integral measurement:
\begin{equation} \label{overdetermination}
\int_\Omega u(t, x, y) \omega(y) dy = \psi(t, x), \quad t \in [0, T], \quad x \in G,
\end{equation}
where $\omega$ and \( \psi \) are known functions, the conditions on which we will determine later.

This paper is organized into seven sections. The next section is auxiliary, where the known results of A. Alikhanov are presented. In section 3 we recall some properties of the eigenvalues of the Laplace operator with the Dirichlet condition. In section 4 a weak solution of the inverse problem under consideration is defined and the main result of the study is presented. The solution of the problem is sought in the form of a series in eigenfunctions with unknown coefficients and in section 5 an infinite system of integro-differential equations is derived to determine these coefficients and a priori estimates are established. Sections 6 and 7 are devoted to the proof of the main result. Finally, section 8 provides a conclusion.

\section{Preliminaries}
In this section, we recall the definition of Mittag-Leffler functions and remind the well-known results of A. Alikhanov, which will be used in further discussions.

The solution of the subdiffusion equations is expressed through the two-parameter Mittag-Leffler functions $E_{\rho,\mu}(z)$ which are determined by the formula
$$
E_{\rho,\mu}(z)= \sum\limits_{k=0}^\infty \frac{z^k}{\Gamma(\rho
k+\mu)},\quad \rho>0, \quad \mu, z\in \mathbb{C}.
$$
If \(\mu = 1\), then we obtain the classical Mittag-Leffler function and is denoted by \(E_{\rho}(z) = E_{\rho,1}(z)\). Obviously, since the Mittag-Leffler functions are an entire function, there exists a constant $M$ such that
\begin{equation}\label{E}
   E_{\rho,\mu}(z) \leq M,  \quad \mu>0, \quad z\in [a,b]\subset \mathbb{R}_+. 
\end{equation}

The following connection between the fractional integral and derivative follows directly from the definitions and equality for $v\in L_1(0,T)$: $J^\alpha_t(J^\beta_t v(t))= J^\beta_t(J^\alpha_t v(t))=J^{\alpha+\beta}_t v(t))$, $\alpha, \beta >0$ (see \cite{Dzh66}, p. 567).
\begin{лем}\label{J} Let $0<\alpha<1$ and $v(t)$ be absolutely continuous on $[0,T]$: $v\in AC\,[0,T]$. Then
\begin{equation}\label{JD}
    J_t^\alpha D^\alpha_t v(t)= v(t)- v(0).
\end{equation}
    
\end{лем}
Proof see in \cite{Dzh66}, p. 570.

Let us now quote the following two statements, in a form convenient for us, from  Alikhanov \cite{Alixan}.
\begin{лем}\label{Alixan1}
Let $0<\alpha<1$ and  $w\in AC\, ([0,T]: L_2(G))$. Then
\[
\frac{1}{2} D_t^\alpha \| w(t,x) \|_{L_2(G)}^2 \leq \int_G w(t,x) \, D_t^\alpha w(t,x) \, dx.
\]
\end{лем}

\begin{лем}\label{Alixan2L}
Let \( y(t)\in AC\, [0,T] \) be a positive function and for all \( t \in (0,T] \), the following inequality holds:
\[
D_t^\alpha y(t) \leq c_1 y(t) + c_2(t), \quad 0<\alpha\leq 1,
\]
for almost all $t\in [0,T]$, where \( c_1>0 \) and \( c_2(t) \) is an integrable nonnegative function on $[0,T]$. Then
\begin{equation}\label{Alixan2}
y(t)\leq y(0) \,E_\alpha(c_1 t^\alpha) + \Gamma(\alpha) \,E_{\alpha, \alpha}(c_1 t^\alpha) \, J_t^\alpha c_2(t).
    \end{equation}
\end{лем}

Note also that if $v(t)$ is a positive bounded function on $[0, T]$ and $\alpha\in (0,1]$, then the following relations hold:
\begin{equation}\label{integral}
\frac{1}{T^{1-\alpha}}\int_0^tv(s) ds \leq \int_0^t\frac{v(s)}{t^{1-\alpha}} ds\leq \int_0^t\frac{v(s)}{(t-s)^{1-\alpha}}ds=\Gamma(\alpha) J_t^\alpha v(t)\leq T^\alpha \max\limits_{t\in [0,T]} v(t).
    \end{equation}
    
Next, we introduce some function spaces. Let \( B \) be a Banach space. We denote by \( L_\infty(0, T; B) \) the space of functions that are essentially bounded on \( (0, T) \) and take values in \( B \). The space \( L_1(0, T; B) \) is defined similarly. Let \( W_2^k(D) \) denote a classical Sobolev space. Then, the symbol \( \dot{W}_2^k(D) \) represents the closure of the set \( C_0^\infty(\Omega) \) with respect to the norm of \( W_2^k(\Omega) \). 
    
    \section{Eigenvalues of the Laplace Operator with Dirichlet Boundary Conditions}
\label{sec:eigenvalues}

Let \(\Omega \subset \mathbb{R}^n\) be a bounded domain with a smooth boundary \(\partial \Omega\). Consider the Laplace operator \(-\Delta\), equipped with Dirichlet boundary conditions, i.e., \(u = 0\) on \(\partial \Omega\). The eigenvalue problem for the Laplace operator is given by:
\begin{equation}\label{spectral_problem}
-\Delta v = \lambda v \quad \text{in } \Omega, \quad v = 0 \quad \text{on } \partial \Omega.
\end{equation}
This operator's spectrum consists of a sequence of positive eigenvalues \(0 < \lambda_1 \leq \lambda_2 \leq \lambda_3 \leq \dots\), counted with multiplicity, such that \(\lambda_k \to +\infty\) as \(k \to \infty\). Each eigenvalue \(\lambda_k\) corresponds to an eigenfunction \(v_k \), and the eigenfunctions \(\{v_k\}\) form an orthonormal basis in \(L^2(\Omega)\).

The asymptotic behavior of the eigenvalues is governed by Weyl's law, which provides insight into their growth rate. For large \(k\), the eigenvalues satisfy:
\[
\lambda_k \sim C_n \left( \frac{k}{|\Omega|} \right)^{2/n},
\]
where \(|\Omega|\) denotes the Lebesgue measure of \(\Omega\), and \(C_n = 4\pi^2 (B_n)^{-2/n}\), with \(B_n = \pi^{n/2} / \Gamma(n/2 + 1)\) being the volume of the unit ball in \(\mathbb{R}^n\). This implies that \(\lambda_k\) grows approximately as \(k^{2/n}\). Weyl's law is a cornerstone of spectral theory and has been extensively studied (see, e.g., \cite{ReedSimon1978, Weyl1912}).

Consider the series involving the eigenvalues \(\lambda_k\) of the Laplace operator with Dirichlet boundary conditions:
\[
    \sum_{k=1}^\infty \lambda_k^{-a},
\]
where \(a > 0\) is a parameter. The convergence of this series depends on the asymptotic behavior of \(\lambda_k\). Since \(\lambda_k \sim k^{2/n}\), then for large \(k\),
\[
\lambda_k^{-a} \sim (k^{2/n})^{-a} = k^{-2a/n}.
\]
The series \(\sum_{k=1}^\infty \lambda_k^{-a}\) is therefore comparable to the p-series \(\sum_{k=1}^\infty k^{-2a/n}\). A p-series \(\sum k^{-p}\) converges if and only if \(p > 1\). Applying this criterion, we require:
\[
\frac{2a}{n} > 1 \implies a > \frac{n}{2}.
\]
Thus, the series \(\sum_{k=1}^\infty \lambda_k^{-a}\) converges if \(a > n/2\) and diverges if \(a \leq n/2\). For any positive $\varepsilon >0$ we denote
\begin{equation}\label{Cvarepsilon}
C_\varepsilon:= \sum_{j=1}^\infty \lambda_j^{-\frac{n}{2} - \varepsilon}.
\end{equation}

Let us derive an estimate for the \( L^2 \)-norm of the gradient of the eigenfunctions \( v_k \) of the Laplace operator with Dirichlet boundary conditions. Let \(\lambda_k\) be the eigenvalues and \( u_k \) the corresponding eigenfunctions satisfying:
\[
-\Delta v_k = \lambda_k v_k \quad \text{in } \Omega, \quad v_k = 0 \quad \text{on } \partial \Omega,
\]
with \(\| v_k \|_{L^2(\Omega)} = 1\). Multiply this equation by \( v_k \) and integrate over \(\Omega\):
\[
-\int_\Omega v_k \Delta v_k \, dy = \lambda_k \int_\Omega v_k^2 \, dy = \lambda_k,
\]
since \(\| v_k \|_{L^2(\Omega)} = 1\). Using Green's identity and the Dirichlet boundary condition (\( v_k = 0 \) on \(\partial \Omega\)), the left-hand side becomes:
\[
-\int_\Omega v_k \Delta v_k \, dy = \int_\Omega |\nabla v_k|^2 \, dy,
\]
as the boundary term vanishes. Thus:
\begin{equation}\label{gradient}
    \int_\Omega |\nabla v_k|^2 \, dy = \lambda_k.
\end{equation}

Again, by Green's identity and since \(\| v_k \|_{L^2(\Omega)} = 1\), we have
\begin{equation} \label{nabla}
\int_\Omega \nabla v_k(y) \nabla v_j(y)  \, dy =- \int_\Omega \Delta v_k(y) v_j(y)  \, dy = \lambda_k \int_\Omega v_k(y) v_j(y)  \, dy =\begin{cases}
		\lambda_k, \,\, \text{if} \,\, j=k;\\
        
         0,\,\, \text{if} \,\, j\neq k.
	\end{cases}
\end{equation} 

Let $A$ be a self-adjoint extension of the Laplace operator with domain $D(-\Delta) = \{v\in C^2(\Omega): v(y)=0, \,\, y\in \partial \Omega\}$. Then $A$ is a positive operator with domain $W_2^2(\Omega) \cap \dot{W}_2^1(\Omega) $ (see, for example, \cite{Berezans’kii}, Chapte 2). Therefore, we can define fractional powers of this operator using the von Neumann theorem. Namely, let $\tau$ be an arbitrary non-negative number. Then
\[
A^\tau p(y) = \sum_{k=1}^\infty \lambda_k^\tau\,  p_k\,  v_k(y),
\]
and domain of the operator $A^\tau$ has the form
\[
D(A^\tau) = \{p\in L_2(\Omega):  \sum_{k=1}^\infty \lambda_k^{2\tau}\, |p_k|^2 < \infty\}.
\]
Here after, the symbol \( p_k \) denotes the Fourier coefficients of a function \( p(y) \in L_2(\Omega) \) with respect to the system \( \{v_k(y)\} \), defined as a scalar product in $L_2(\Omega)$: $p_k=(p, v_k)$.

Note, that if $\tau_1> \tau_2$ then $D(A^{\tau_1})\subset D(A^{\tau_2}) $. $D(A^\tau)$ is a linear space and we can introduce a scalar product in it using the formula 
\[
(p, q)_{\tau} = \sum_{k=1}^\infty \lambda_k^{2\tau}\,  p_k\, q_k\, =\, (A^\tau p, A^\tau q),
\]
which defines the following norm
\[
||p||^2_{\tau} = \sum_{k=1}^\infty \lambda_k^{2\tau}\,  |p_k|^2\, = \, \int_\Omega |A^\tau p(y)|^2 \, dy.
\]
In many cases in this paper, the function $p$ also depends on $x\in G$ and on $t\in [0, T]$, i.e. $p= p(t,x, y)$. In such cases, we introduce the notation
\[
||p(t, \cdot, \cdot)||_{\tau, G}^2 :=\, \int_G \int_\Omega |A^\tau p(t,x,y)|^2 dy\, dx = \int_G \sum_{k=1}^\infty \lambda_k^{2\tau}\,  |p_k (t, x)|^2\, dx. 
\]

The question naturally arises: how to check whether a given function $p$ belongs to $D(A^\tau)$?	This question has been thoroughly studied, for example, in the work of D. Fujiwara \cite{Fujiwara} (for high-order operators, see the work of Sh. Alimov \cite{Alimov}). If the order $\tau$ is an integer, then the answer to this question is much simpler. For example, let us define the domain of the operator \((-\Delta)^2\). First note that $D(\Delta^{2})\subset D(-\Delta) $. Therefore, if \(p\in D(\Delta^2)\), then $p\in W_2^2(\Omega) \cap \dot{W}_2^1(\Omega) $ and $\Delta p \in W_2^2(\Omega) \cap \dot{W}_2^1(\Omega)$. Or in other words,
\[
D((-\Delta)^2) = W_2^4(\Omega) \cap \{ u \in \dot{W}_2^1(\Omega) \mid \Delta u \in \dot{W}_2^1(\Omega) \}.
\]
If the equality to zero on the boundary of the domain $\Omega$ is understood as the equality to zero of the traces, then we can write
\begin{equation}\label{Delta^2}
    D((-\Delta)^2) = \{ u \in W_2^4(\Omega) \mid u = 0, \Delta u = 0 \text{ on } \partial \Omega \}.
\end{equation}

\section{Definition of the Weak Solution and Formulation of the Main Result}First, we find a formal representation for the unknown function $h$ using the overdetermination condition \eqref{overdetermination}. To do this, we represent the solution to the forward problem (\ref{1}) as a formal series in terms of the eigenfunctions \(\{ v_k (y)\} \) of the spectral problem (\ref{spectral_problem}):
\begin{equation}\label{6}
    u(t, x, y) = \sum_{k=1}^\infty u_k(t, x) v_k (y),
\end{equation}
where \( u_k(t, x) \) are the unknown functions.

In order to use the overdetermination condition (\ref{overdetermination}), we multiply the equation in (\ref{1}) by the function $\omega$ and integrate over the domain $\Omega$. Then
\begin{equation*}
    D_t^{\alpha} \psi - \Delta_x \psi - (\Delta_y u, \omega) = (f(t, x, \cdot), \omega) \,h(t, x) + (g(t, x, \cdot), \omega).
\end{equation*}
Let \( (f(t, x, \cdot), \omega) \neq 0 \) for all \( (t, x) \in [0, T] \times G \). Then the desired representation for \( h(t, x) \) can be written as:
\begin{equation}\label{7}
    h(t, x) = \frac{D_t^{\alpha} \psi(t, x) - \Delta_x \psi(t, x)-(g(t, x, \cdot), \omega) - (\Delta_y u, \omega)}{(f(t, x, \cdot), \omega)}.
\end{equation}

Next, we will investigate the following weak formulation of the inverse problem \eqref{1}–\eqref{overdetermination}.

\begin{опр}\label{opr1} 
Find a pair of functions \( \{ u(t, x, y), h(t, x) \} \), where \( h(t, x) \) has the form \eqref{7}, and the function \( u(t, x, y) \) satisfies the following conditions:
\begin{enumerate}
    \item \( u \in L_\infty(0, T; L_2(D)) \), \,\, \( u \in L_1(0, T; \dot{W}_2^1(D)) \);
    \item \( D_t^\alpha u \in L_1(0, T; L_2(D)) \);
    \item \( u(0, x, y) = \varphi(x, y) \) a.e.\ in \( D \);
    \item For any \( v \in \dot{W}_2^1(D) \) and almost every \( t \in (0, T] \), the following equality holds:
\end{enumerate}
\begin{equation}\label{equation}
    \int_D D_t^\alpha u v \, dx \, dy + \int_D (\nabla_x u \nabla_x v + \nabla_y u \nabla_y v) \, dx \, dy = \int_D f h v \, dx \, dy + \int_D g v \, dx \, dy.
\end{equation}
\end{опр}

Let $M_\alpha = M T^\alpha$, where $M$ is from (\ref{E}) and recall the constant $C_\varepsilon$ is defined in (\ref{Cvarepsilon}). Next, it is convenient for us to introduce a notation for some positive number $\varepsilon$: 
\begin{equation}\label{tau}
2\tau:= 2\tau(\varepsilon) = \frac{n}{2} +1 + \varepsilon.
\end{equation}

We now state the main result of the paper.

\begin{thm}\label{theorem 1}

Let $\varepsilon$ be any positive number and suppose the following conditions hold:
\begin{enumerate}
    \item \( f(t, x, y) \in C(Q) \), \( (f(t, x, \cdot), \omega) \neq 0 \), and \( f_0 = \max_{t, x} \left| \frac{1}{(f(t, x, \cdot), \omega)} \right| \);
    \item \( g(t, x, y) \in C(Q) \), and \( g_0 = \max_{t, x} |(g(t, x, \cdot), \omega)| \);
    \item \( D_t^\alpha \psi, \,\Delta \psi \in C([0, T] \times G) \), and \( \psi_0 = \max_{t, x} \left( |D_t^\alpha \psi| + |\Delta\psi| \right) \);
    \item   \( M_\alpha C_\varepsilon\, f_0^2\, ||\nabla \omega||_{L_2(\Omega)}^2 \max_{t,x} ||f(t, x, \cdot)||^2_{\tau} \leq 1\), and \(f^*=\max_t||f(t, \cdot, \cdot)||_{\tau,\, G}^2 < \infty\);
    \item \( g^*=\max_t||g(t, \cdot, \cdot)||_{\tau,\, G}^2 < \infty \);
     \item \( \|\nabla_x \varphi \|^2_{L_2(D)} < \infty, \) and \(
\varphi^*= ||\varphi||_{\tau,\, G}^2 < \infty\);
     \item \(\omega \in \dot{W}^1_2(\Omega)\) and \(\int_\Omega \varphi(x,y) \omega(y) dy= \psi(0, x)\).
\end{enumerate}
Then, the inverse problem has a unique weak solution. Moreover, the following estimates hold:
\begin{equation}\label{uestimate1}
\max_{[0,T]}\| u(t, \cdot, \cdot) \|^2_{\tau, G}\leq 2 A_0.
\end{equation}
 \begin{equation}\label{u_xu_y}
||\nabla_x u ||_{L_2(Q)}^2 +||\nabla_y u ||_{L_2(Q)}^2 \leq \frac{\Gamma(\alpha)T^{1-\alpha}}{2}||\varphi||^2_{L_2(D)}+ 3 A_0 T +\frac{T}{2} \, A_1 
\end{equation}
\[
+T\, A_0 C_\varepsilon\, f_0^2\, ||\nabla \omega||_{L_2(\Omega)}^2 \max_{x\in G} ||f(t,x,\cdot)||_{L_2(\Omega)}^2, 
\]
\begin{equation}\label{DalphaInt}
    || D_t^\alpha u||_{L_2(Q)}^2 \leq \frac{\Gamma(\alpha)T^{1-\alpha}}{2} \left(\|\nabla_x \varphi \|^2_{L_2(D)} +  \|\nabla_y \varphi \|^2_{L_2(D)}\right)
\end{equation}
\[
+T\max_t\left[f_0 ^2(\psi_0 + g_0)^2 \|f(t, \cdot, \cdot)\|^2_{L_2(D)}
+ \|g(t, \cdot, \cdot)\|^2_{L_2(D)}\right]
\]
\[
+2\, A_0 C_\varepsilon\, f_0^2\, ||\nabla \omega||_{L_2(\Omega)}^2 \max_{x\in G} ||f(t,x,\cdot)||_{L_2(\Omega)}^2,
\]

\[
   ||h(t,\cdot)||_{L_2(G)}^2 \leq 4 f_0^2 \left( ||D_t^{\alpha} \psi(t, x)||_{L_2(G)}^2 + ||\Delta_x \psi(t, x)||_{L_2(G)}^2 + ||(g(t, x, \cdot), \omega)||_{L_2(G)}^2\right)
   \]
   \[
   + 8 f_0^2\, A_0\, C_\varepsilon ||\nabla \omega||_{L_2(\Omega)}^2.
\]
Here  
\[
A_1= f_0^2 (\psi_0 + g_0)^2 \max_{t}\|f(t, \cdot, \cdot)\|_{L_2(D)}^2 + \max_{t}\|g(t, \cdot, \cdot)\|_{L_2(D)}^2,
\]
and 
\[
A_0 = M \varphi^* + M_\alpha f_0^2(\psi_0 + g_0)^2 f^*+ M_\alpha g^*. 
\]

\end{thm}

\begin{зам}\label{zam.1} It is known that for $\alpha>n/4$ the eigenfunction expansion of any function $f\in D(A^\alpha)$ converges absolutely and uniformly in the closed domain $\overline{\Omega}$ (see \cite{Ilyin1958}). And the functions $\varphi$, $f$ and $g$ belong to the class $D(A^\tau)$, with $\tau> \frac{n+2}{4}$, which naturally guarantees the indicated convergence.
\end{зам}

\section{A priori estimates}

Let us write the function $h$ in a form convenient for us. First, taking into account the property of the function $\omega$ (see condition 7 of the theorem), we get $(\Delta_y u, \omega) = (\nabla_y u, \nabla \omega)$. Then, applying the formal form \eqref{6} of the function $u$, we will have
\begin{equation}\label{h_final}
    h(t, x) = \frac{D_t^{\alpha} \psi(t, x) - \Delta_x \psi(t, x)-(g(t, x, \cdot), \omega) - \sum_{j=1}^\infty (\nabla v_j, \nabla \omega) u_j(t,x) }{(f(t, x, \cdot), \omega)}.
\end{equation}

We decompose the functions \( f(t, x, y) \), \( g(t, x, y) \), and \( \varphi(x, y) \) into Fourier series (see Remark \ref{zam.1}) and denote their corresponding Fourier coefficients by \( f_k(t, x) \), \( g_k(t, x) \), and \( \varphi_k(x) \). Substitute these series and the representation \eqref{6} for \( u(t, x, y) \) into equation \eqref{equation}, with the test function \( v(x, y) \) replaced by \( w(x) v_k (y) \), where \( w \in \dot{W}_2^1(G) \). Next note that according to (\ref{nabla}) one has
\begin{equation}\label{gradientlambda}
\int_D\nabla_y u \nabla_y v dx dy= \sum_{j=1}^\infty \int_G u_j(t,x) w(x) dx \int_\Omega \nabla v_j \nabla v_k dy  = \lambda_k \int_G u_k w \, dx.
\end{equation}
Thus we obtain the following equation to determine the unknown coefficients \( u_k(t, x) \) (see \eqref{6}):
\begin{equation}\label{81}
\int_G D_t^\alpha u_k w \, dx + \int_G \nabla u_k \nabla w \, dx + \lambda_k \int_G u_k w \, dx = \int_G (g_k + f_k h) w \, dx,
\end{equation}
but here \( h \) is defined using all \( u_j(t, x) \), \( j = 1, 2, \dots \) (see \eqref{h_final}). This fact creates a certain problem and in order to solve it we use the method of successive approximations, i.e. using recurrence relations we construct a sequence \( \{ u_k^i \} \), \( i = 1, 2, \dots \), and then we prove that \( \{ u_k^i \} \to u_k \) as \( i \to \infty \) in the appropriate norm.

The corresponding recurrence relations for all \( w \in \dot{W}_2^1(G) \), \( k \geq 1 \) and \( i \geq 1 \) have the form:
\begin{equation}\label{9}
\int_G D_t^\alpha u_k^i w \, dx + \int_G \nabla u_k^i \nabla w \, dx + \lambda_k \int_G u_k^i w \, dx 
\end{equation}
\[
= \int_G M_k(t, x) w \, dx - \int_G \frac{f_k(t, x)}{(f(t, x, \cdot), \omega)}  \sum_{j=1}^\infty (\nabla v_j, \nabla \omega) u^{i-1}_j(t,x) w \, dx,
\]
with the initial conditions:
\begin{equation}\label{10}
u_k^i(0, x) = \varphi_k(x),
\end{equation}
where
\[
M_k(t, x) = f_k(t, x) \left[ \frac{D_t^\alpha \psi(t, x) - \Delta_x \psi(t, x) - (g(t, x, \cdot), \omega)}{(f(t, x, \cdot), \omega)} \right] + g_k(t, x).
\]

If the right-hand side of equation \eqref{9} is known, then the strong formulation of the problem \eqref{9}–\eqref{10} is well studied (see, e.g., \cite{Sakamoto}, \cite{AMux2}).  For instance, in \cite{Sakamoto}, it is established that if the right-hand side (for all \( t \in (0, T] \)) and the initial function \( \varphi_k(x) \) belong to \( L_2(G) \) (which holds in our case), then there exists a unique strong solution to the problem satisfying \( u_k^i(t, x) \in L_\infty(0, T; \dot{W}_2^2(G)) \) and \( D_t^\alpha u_k^i(t, x) \in L_\infty(0, T; L_2(G)) \). Clearly, such a strong solution is also a weak solution.

Here’s a refined version of your text:

We analyze the problem \eqref{9}–\eqref{10} as follows. We set \( u_k^0 = 0 \) as the initial approximation for all \( k = 1, 2, \dots \). Then, we solve the initial-boundary value problem \eqref{9}–\eqref{10} to obtain the sequence \( \{ u_k^1 \} \), \( k = 1, 2, \dots \). Subsequently, we construct the sequence \( \{ u_k^i \} \) iteratively. Prior to this, assuming that all \( \{ u_k^i \} \) have been constructed, we first derive a priori estimates.

\begin{лем}\label{u_k} Let $\varepsilon$ be any possitive number. If
\begin{equation}\label{fk}
M_\alpha C_\varepsilon\, f_0^2\, ||\nabla \omega||_{L_2(\Omega)}^2 \max_{t,x} ||f(t, x, \cdot)||^2_{\tau(\varepsilon)} \leq 1.
\end{equation}
then
\begin{equation}\label{181}
\max_{[0,T]}\sum_{k=1}^\infty \lambda_k^{2 \tau(\varepsilon)} ||u_k^{i}||_{L_2(G)}^2  \leq 2 A_0,
\end{equation}
where 
\[
A_0 = M \varphi^* + M_\alpha f_0^2(\psi_0 + g_0)^2 f^*+ M_\alpha g^*, \,\,
\varphi^*= ||\varphi||_{\tau(\varepsilon), G}^2, 
\]
and
\[
f^*=\max_{[0,T]}||f(t, \cdot, \cdot)||_{\tau(\varepsilon), G}^2 \,\,, \quad g^*=\max_{[0,T]}||f(t, \cdot, \cdot)||_{\tau(\varepsilon), G}^2.
\]

\end{лем}

\begin{proof} Substitute \( w = u_k^i \) into equation \eqref{9}, to obtain:
\begin{equation}\label{82}
\int_G D_t^\alpha u_k^i u_k^i \, dx + \| \nabla u_k^i \|_{L_2(G)}^2 + \lambda_k \| u_k^i \|_{L_2(G)}^2
\end{equation}
\[
 = \int_G M_k(t, x) u_k^i \, dx - \int_G u_k^i\frac{f_k(t, x)}{(f(t, x, \cdot), \omega)}  \sum_{j=1}^\infty (\nabla v_j, \nabla \omega) u^{i-1}_j(t,x)  \, dx,
\]
Apply for the integral on the left-hand side Alikhanov’s estimate (see Lemma \ref{Alixan1}). Then
\begin{equation}\label{12}
\frac{1}{2} D_t^\alpha \| u_k^i \|_{L_2(G)}^2 + \| \nabla u_k^i \|_{L_2(G)}^2 + \lambda_k \| u_k^i \|_{L_2(G)}^2 
\end{equation}
\[
\leq \left| \int_0^1 M_k(t, x) u_k^i \, dx \right| + \left| \int_G u_k^i\frac{f_k(t, x)}{(f(t, x, \cdot), \omega)}  \sum_{j=1}^\infty (\nabla v_j, \nabla \omega) u^{i-1}_j(t,x)  \, dx \right| := I_1(t) + I_2(t).
\]
For the right-hand side, using the notation and conditions of Theorem \ref{theorem 1}, we have:
\[
I_1(t) \leq f_0 (\psi_0 + g_0) \int_G |f_k(t, x) u_k^i(t, x)| \, dx + \int_G |g_k(t, x) u_k^i(t, x)| \, dx.
\]
Apply the inequality \( 2ab \leq a^2 + b^2 \), to get
\begin{equation}\label{13}
I_1(t) \leq \frac{1}{2}\left[f_0^2 (\psi_0 + g_0)^2 \| f_k \|_{L_2(G)}^2 + \| g_k \|_{L_2(G)}^2\right]+  \| u_k^i \|_{L_2(G)}^2.
\end{equation}
For $I_2$ it is not hard to see, that
\[
I_2(t)\leq \frac{1}{2} \| u_k^i \|_{L_2(G)}^2 +  \frac{1}{2}f_0^2 \int_G |f_k(t, x)|^2 \left| \sum_{j=1}^\infty (\nabla v_j, \nabla \omega) u^{i-1}_j(t,x) \right|^2 \, dx.
\]
Let us estimate the sum separately. Due to (\ref{gradient}) and the Cauchy–Bunyakovsky inequality we have
\[
I :=  \left| \sum_{j=1}^\infty (\nabla v_j, \nabla \omega) u^{i-1}_j(t,x) \right|^2 \,  \leq  \left[ \sum_{j=1}^\infty ||\nabla v_j||_{L_2(\Omega)}||\nabla \omega||_{L_2(\Omega)} |u^{i-1}_j(t,x)| \right]^2 
\]
\[
= ||\nabla \omega||_{L_2(\Omega)}^2  \left[ \sum_{j=1}^\infty\lambda_j^\frac{1}{2} |u^{i-1}_j(t,x)| \right]^2.
\]
For any \( \varepsilon > 0 \), we express \( \lambda_j^{\frac{1}{2}} = \lambda_j^{\frac{1}{2} + \frac{n + 2\varepsilon}{4}} \lambda_j^{-\frac{n + 2\varepsilon}{4}} \) and use the Cauchy–Schwarz inequality on the sum. Then, refer to (\ref{Cvarepsilon}):
\begin{equation}\label{C}
I \leq  C_\varepsilon ||\nabla \omega||_{L_2(\Omega)}^2 \sum_{j=1}^\infty \lambda_j^{\frac{n}{2} + 1+\varepsilon} |u_j^{n-1}|^2. 
\end{equation}

By сombining estimates of $I_1,\, I_2$ and (\ref{C}), the expression in (27) can be reformulated as: \begin{equation}\label{15}
\frac{1}{2} D_t^\alpha \| u_k^i \|_{L_2(G)}^2 + \| \nabla u_k^i \|_{L_2(G)}^2 + \lambda_k \| u_k^i \|_{L_2(G)}^2\leq \frac{3}{2} \| u_k^i \|_{L_2(G)}^2   
\end{equation}
\[
 +  \frac{1}{2}  \left[f_0^2 (\psi_0 + g_0)^2 \|f_k(t, \cdot)\|_{L_2(G)}^2 + \|g_k(t, \cdot)\|_{L_2(G)}^2\right]
\]
\[
 + \frac{ C_\varepsilon\, f_0^2\, ||\nabla \omega||_{L_2(\Omega)}^2}{2} \int_G |f_k(t, x)|^2 \sum_{j=1}^\infty \lambda_j^{2 \tau} |u_j^{i-1}|^2 \, dx.
\]
If we omit the final two terms on the left-hand side, then
\[
D_t^\alpha \| u_k^i \|_{L_2(G)}^2 \leq 3\| u_k^i \|_{L_2(G)}^2 +c^k_2(t),
\]
where $ c_2^k(t)= c^k_{2,1}(t)+ c^k_{2, 2}(t)$ and
\[
 c^k_{2,1}(t)= f_0^2 (\psi_0 + g_0)^2 \|f_k(t, \cdot)\|_{L_2(G)}^2 + \|g_k(t, \cdot)\|_{L_2(G)}^2,
 \]
 \[
c^k_{2, 2}(t)= C_\varepsilon\, f_0^2\, ||\nabla \omega||_{L_2(\Omega)}^2 \int_G |f_k(t, x)|^2 \sum_{j=1}^\infty \lambda_j^{2 \tau} |u_j^{i-1}|^2 \, dx.
\]
By Lemma \ref{Alixan2L} we obtain
\[
\| u_k^i \|_{L_2(G)}^2\leq \| \varphi_k \|_{L_2(G)}^2 \,E_\alpha(3 t^\alpha) + \Gamma(\alpha) \,E_{\alpha, \alpha}(3 t^\alpha) \, J_t^\alpha c^k_2(t).
\]
Relying on the boundedness of the Mittag-Leffler function (refer to (\ref{E})) and the estimate (\ref{integral}), we deduce that
\[
\| u_k^i(t, \cdot) \|_{L_2(G)}^2\leq M\| \varphi_k \|_{L_2(G)}^2  + MT^\alpha \left(\max_{[0,T]}c^k_{2,1}(t) +\max_{[0,T]}c^k_{2,2}(t)\right).
\]
By multiplying this inequality by $\lambda_j^{2 \tau}$ and summing over $k$ from 1 to $\infty$, it follows that
\begin{align*}
&\max_{[0,T]}\sum_{k=1}^{\infty} \lambda_k^{2 \tau} \| u_k^i (t, \cdot)\|^2_{L_2(G)} \leq  M \varphi^* + M_\alpha f_0^2(\psi_0 + g_0)^2 f^* + M_\alpha g^*\\
&+ M_\alpha C_\varepsilon\, f_0^2\, ||\nabla \omega||_{L_2(\Omega)}^2\max_{[0,T]}\int_G \sum_{k=1}^{\infty} \lambda_k^{2 \tau} |f_k(t, x)|^2 \sum_{j=1}^{\infty} \lambda_j^{2 \tau} | u_j^{i-1}(t, x)|^2 dx.
\end{align*}
Employing condition (\ref{fk}) and the notation from Lemma \ref{u_k}, the final estimate can be reformulated as a recurrence estimate:
\begin{equation}\label{17}
\max_{[0,T]}\sum_{k=1}^{\infty} \lambda_k^{2 \tau}\| u_k^i \|_{L_2(G)}^2 \leq A_0 + \frac{1}{2} \max_{[0,T]}\sum_{k=1}^{\infty} \lambda_k^{2 \tau} \| u_k^{i-1} \|_{L_2(G)}^2, \quad i = 1, 2, \dots.
\end{equation}

As stated previously, we set \( u_k^0 = 0 \) as the initial approximation for all \( k \geq 1 \). Consequently, for \( u_k^1(t, x) \), \( k = 1, 2, \dots \), using \eqref{17}, we derive:
\begin{equation}\label{u1}
\max_{[0,T]}\sum_{k=1}^\infty \lambda_k^{2 \tau} \| u_k^1 \|_{L_2(G)}^2 \leq A_0.
\end{equation}
Subsequently, we insert the functions \( \{ u_k^1 \}_{k=1}^\infty \) into the problem \eqref{9}–\eqref{10} to uniquely define the functions \( \{ u_k^2 \}_{k=1}^\infty \). For these functions, we obtain the following estimate from \eqref{17}:
\[
\max_{[0,T]}\sum_{k=1}^\infty \lambda_k^{2 \tau} \| u_k^2 \|_{L_2(G)}^2 \leq A_0 + \frac{1}{2} A_0.
\]
Proceeding with this approach, we ultimately derive:
\begin{equation}\label{18}
\max_{[0,T]}\sum_{k=1}^\infty \lambda_k^{2 \tau} \| u_k^n \|_{L_2(G)}^2 = A_0 \sum_{j=1}^n \left( \frac{1}{2} \right)^j = 2 \left( 1 - \frac{1}{2^n} \right) A_0 \leq 2 A_0.
\end{equation}
\end{proof}

\begin{лем}The following estimate is valid
\[
\sum\limits_{k=1}^\infty\int\limits_0^t\int\limits_G\left[| \nabla\, u_k^i |^2 + \lambda_k | u_k^i |^2\right]dx d\tau \leq \frac{\Gamma(\alpha)T^{1-\alpha}}{2}||\varphi||^2_{L_2(D)}+ 3 A_0 T +\frac{T}{2} \, A_1
\]
\begin{equation}\label{u_x}
+T\, A_0 C_\varepsilon\, f_0^2\, ||\nabla \omega||_{L_2(\Omega)}^2 \max_{x\in G} ||f(t,x,\cdot)||_{L_2(\Omega)}^2, \quad t\in [0,T]. 
\end{equation}

\end{лем}

\begin{proof}

Applying the operator $J_t^\alpha$ to both sides of inequality (\ref{15}) and utilizing Lemma \ref{J}, we obtain:
\begin{equation}\label{151}
\| u_k^i \|_{L_2(G)}^2 + 2 J_t^\alpha\| \nabla u_k^i \|_{L_2(G)}^2 + 2 \lambda_k J_t^\alpha\| u_k^i \|_{L_2(G)}^2 
\end{equation}
\[
\leq \| \varphi_k \|_{L_2(G)}^2+3 J_t^\alpha \| u_k^i \|_{L_2(G)}^2 +  J_t^\alpha c_{2,1}^k(t) +J_t^\alpha c_{2,2}^k(t).
\]
Initially, we disregard the first term on the left-hand side of the inequality. Then, by summing inequality \eqref{151} over \( k \) from 1 to \( \infty \) and applying the estimates (\ref{integral}), we derive:\begin{equation}\label{151}
\frac{2}{\Gamma(\alpha)T^{1-\alpha}} \sum\limits_{k=1}^\infty\int\limits_0^t\int\limits_G\left[| \nabla u_k^i|^2 +  \lambda_k | u_k^i |^2\right]dxd\tau 
\end{equation}
\[
\leq \sum\limits_{k=1}^\infty\| \varphi_k \|_{L_2(G)}^2+\frac{3 T^\alpha}{\Gamma(\alpha)}\max_{[0,T]} \sum\limits_{k=1}^\infty\| u_k^i \|_{L_2(G)}^2 +  \frac{ T^\alpha}{\Gamma(\alpha)}\max_{[0,T]} \sum\limits_{k=1}^\infty \left(c_{2,1}^k(t) + c_{2,2}^k(t)\right).
\]
Apply Parseval's equality to obtain
\[
\max_{[0,T]}\sum\limits_{k=1}^\infty c_{2,1}^k(t)=f_0^2 (\psi_0 + g_0)^2 \max_{[0,T]}\|f(t, \cdot, \cdot)\|_{L_2(D)}^2 + \max_{[0,T]}\|g(t, \cdot, \cdot)\|_{L_2(D)}^2:= A_1.
\]
Similarly,
\[
\sum\limits_{k=1}^\infty c_{2,2}^k(t) \leq C_\varepsilon\, f_0^2\, ||\nabla \omega||_{L_2(\Omega)}^2 \max_{x\in G}||f(t, x, \cdot)||_{L_2(\Omega)}^2\, \sum_{j=1}^\infty \lambda_j^{2 \tau} ||u_j^{i-1}||_{L_2(G)}^2,
\]
or by (\ref{181}),
\[
\sum\limits_{k=1}^\infty c_{2,2}^k(t) \leq 2\, A_0 C_\varepsilon\, f_0^2\, ||\nabla \omega||_{L_2(\Omega)}^2 \max_{x\in G} ||f(t, x, \cdot)||_{L_2(\Omega)}^2\,. 
\]
Thus, inequality (\ref{151}) can be reformulated as (\ref{u_x}).
\end{proof}

\begin{лем} For any $t\in [0,T]$ we have
\begin{equation}\label{Dalpha1}
\sum\limits_{k=1}^\infty\int\limits_0^t\int\limits_G| D_t^\alpha u_k^i |^2dxd\tau \leq \frac{\Gamma(\alpha)T^{1-\alpha}}{2} \left(\| \nabla_x\varphi \|^2_{L_2(D)} +  \| \nabla_y\varphi \|_{L_2(D)}^2\right)
\end{equation}
\[
+T\max_t\left[f_0 ^2(\psi_0 + g_0)^2 \|f(t, \cdot, \cdot)\|^2_{L_2(D)}
+ \|g(t, \cdot, \cdot)\|^2_{L_2(D)}\right]
\]
\[
+2\, A_0 C_\varepsilon\, f_0^2\, ||\nabla \omega||_{L_2(\Omega)}^2 \max_{x\in G} ||f(t, x, \cdot)||_{L_2(\Omega)}^2.
\]
\end{лем}
\begin{proof}
Substitute \( w = D_t^\alpha u_k^i \) into equation \eqref{9}, to get
\begin{equation}\label{82}
\| D_t^\alpha u_k^i \|_{L_2(G)}^2 + \int_G \nabla u_k^i\, D_t^\alpha\, \nabla u_k^i \, dx + \lambda_k \int_G u_k^i\, D_t^\alpha\,u_k^i \, dx 
\end{equation}
\[
= \int_G M_k(t, x) D_t^\alpha u_k^i \, dx - \int_G D_t^\alpha u_k^i \frac{f_k(t, x)}{(f(t, x, \cdot), \omega)}  \sum_{j=1}^\infty (\nabla v_j, \nabla \omega) u^{i-1}_j(t,x) \, dx.
\]
We apply A. Alikhanov's estimate (refer to Lemma \ref{Alixan1}) to the second and third integrals on the left-hand side of the equality, yielding:
\begin{equation}\label{121}
   \| D_t^{\alpha} u_k^i \|^2_{L_2(G)} + \frac{1}{2} D_t^{\alpha} \| \nabla u_k^i \|^2_{L_2(G)} + \lambda_k \frac{1}{2} D_t^{\alpha} \| u_k^i \|^2_{L_2(G)}  
\end{equation}
\[
\leq \left|\int_G M_k(t,x) D_t^{\alpha} u_k^i dx\right| + \left|\int_G D_t^{\alpha} u_k^i \frac{f_k(t, x)}{(f(t, x, \cdot), \omega)}  \sum_{j=1}^\infty (\nabla v_j, \nabla \omega) u^{i-1}_j(t,x)  \, dx\right| := J_1(t) + J_2(t).
\]
Apply inequality $a\, b \leq a^2 + \frac{1}{4} \, b^2$ to $J_1(t)$ to get (see Eq. \eqref{13})
\[
J_1(t)  \leq  f_0^2(\psi_0  + g_0)^2 ||f_k||_{L_2(G)}^2
+  ||g_k||_{L_2(G)}^2 + \frac{1}{2} \| D_t^{\alpha} u_k^n \|^2_{L_2(G)}.
\]
The second summand in the right-hand side of \eqref{121} has the estimate (see Eq. \eqref{C})
\[
J_2(t)\leq \frac{1}{4} \| D_t^{\alpha} u_k^i \|^2_{L_2(G)} + C_\varepsilon \, f_0^2 \,||\nabla \omega||_{L_2(\Omega)}^2 \int_G |f_k(t,x)|^2 \sum_{j=1}^\infty \lambda_j^{2 \tau} |u_j^{i-1}|^2   dx.
\]

Now we apply the operator $J_t^\alpha$ to both parts of inequality (\ref{15}) and use estimates of $J_1(t), \, J_2(t)$ and  Lemma \ref{J}. Then
\[
\frac{1}{4}J_t^\alpha\| D_t^{\alpha} u_k^i \|^2_{L_2(G)}+ \frac{1}{2} \| \nabla u_k^i \|^2_{L_2(G)} + \lambda_k \frac{1}{2} \| u_k^i \|^2_{L_2(G)}
\]
\[
\leq \frac{1}{2} \| \nabla \varphi_k \|^2_{L_2(G)} + \lambda_k \frac{1}{2} \| \varphi_k \|^2_{L_2(G)}+f_0 ^2(\psi_0 + g_0)^2 J_t^\alpha\|f_k\|^2_{L_2(G)}
+  J_t^\alpha\|g_k\|^2_{L_2(G)}
\]
\[
+  C_\varepsilon \, f_0^2 \,||\nabla \omega||_{L_2(\Omega)}^2 J_t^\alpha\, \int_G |f_k(t,x)|^2 \sum_{j=1}^\infty \lambda_j^{2 \tau} |u_j^{i-1}|^2   dx.
.
\]
Initially, we exclude the second and third terms on the left-hand side of the inequality. Then, by summing over \( k \) from 1 to \( \infty \) and utilizing the estimates (\ref{integral}), we derive (\ref{Dalpha1}).
\end{proof}

\section{Convergence}
We aim to demonstrate that the sequences \( u_k^i \), \( \nabla u_k^i \), and \( D_t^\alpha u_k^i \), for \( i=1,2,\dots \), are fundamental with respect to their corresponding norms for all \( k \geq 1 \).

Given that the initial-boundary value problem \eqref{9}–\eqref{10} is linear with respect to \( u_k^i \), by reiterating the previous arguments, we can establish estimates \eqref{18}, \eqref{u_x}, and \eqref{Dalpha} for \( u_k^{i+p} - u_k^i \), where \( p = 1,2,3,\dots \). For instance, consider writing equality (\ref{9}) for \( u_k^{i+1} \) and subtracting (\ref{9}) from it. This  leads  in (noting that here \( M_k(t, x) \equiv 0 \)):
\begin{equation}\label{Dalpha}
\int_G D_t^\alpha (u_k^{i+1}-u_k^i)\, w \, dx + \int_G \nabla(u_k^{i+1}-u_k^i) \,\nabla w \, dx + \lambda_k \int_G (u_k^{i+1}-u_k^i) \,w \, dx 
\end{equation}
\[
=  - \int_G \frac{f_k(t, x)}{(f(t, x, \cdot), \omega)}  \sum_{j=1}^\infty (\nabla v_j, \nabla \omega) u^{i-1}_j(t,x) w \, dx,
\]
for all \( w \in \dot{W}_2^1(0, 1) \), with the initial condition:
\begin{equation}\label{101}
(u_k^{i+1}-u_k^i)(0, x) = 0.
\end{equation}
By applying the same  reasoning as in the proof of \eqref{17}, but substituting \eqref{101} for the initial condition \eqref{10}, we obtain the inequality:
\[
\max_{[0,T]}\sum_{k=1}^{\infty} \lambda_k^{2 \tau}\| u_k^{i+1} -u_k^i \|_{L_2(G)}^2 \leq\frac{1}{2}\max_{[0,T]}\, \sum_{k=1}^{\infty} \lambda_k^{2 \tau}\| u_k^{i} -u_k^{i-1} \|_{L_2(G)}^2.
\]
Since \( u_k^0 = 0 \) for all \( k \geq 1 \), then from \eqref{u1}, we get
\[
\max_{[0,T]}\, \sum_{k=1}^{\infty} \lambda_k^{2 \tau}\| u_k^{2} -u_k^1 \|_{L_2(G)}^2 \leq A_0. 
\]
For any \( n \geq 2 \) the last two estimates imply:
\[
\max_{[0,T]}\, \sum_{k=1}^{\infty} \lambda_k^{2 \tau}\| u_k^{i+1} -u_k^i \|_{L_2(G)}^2 \leq A_0 \left( \frac{1}{2} \right)^{n-2}.
\]
Hence, if we write $u_k^{i+p} - u_k^i=u_k^{i+p} - u_k^{i+p-1}+ u_k^{i+p-1}- u_k^{i+p-2}+\cdots$, then
\begin{equation}\label{fundamental}
\max_{[0,T]}\, \sum_{k=1}^{\infty} \lambda_k^{2 \tau}\| u_k^{i+p} -u_k^i \|_{L_2(G)}^2 \leq A_0 \left( \frac{1}{2} \right)^{i-2} \left( 1 + \frac{1}{2} + \dots + \frac{1}{2^{p}} \right)\leq A_0 \left( \frac{1}{2} \right)^{i-1}.
\end{equation}
This indicates that the sequence \(\{u_k^i\}\) is fundamental in \(L_2(G)\). Therefore, for each \(k \geq 1\), there exist functions \(u_k(t,x) \in L_2(G)\) such that, for every \(0 < t \leq T\), as \(i \to \infty\), the following conditions are satisfied:
\begin{equation}\label{converU}
\begin{cases}
u_k^i(t, x) \to u_k (t, x)\quad \text{in } L_2(G), \quad k=1,2, 
\cdots,\\
\sum_{k=1}^{\infty} \lambda_k^{2 \tau}\| u_k^{i}(t, \cdot) -u_k(t, \cdot) \|_{L_2(G)}^2  \to 0.
\end{cases}
\end{equation}

By substituting $ u_k^{i+p} - u_k^i $ for $ u_k^i $ in the estimate \eqref{151}, it assumes the form (here, one must replicate the reasoning analogous to that provided in the proof of (\ref{fundamental})):
\[
\| u_k^{i+p}-u_k^i \|_{L_2(G)}^2 + 2 J_t^\alpha\| \nabla(u_k^{i+p}-u_k^i ) \|_{L_2(G)}^2 + 2 \lambda_k J_t^\alpha\| u_k^{i+p}-u_k^i  \|_{L_2(G)}^2 
\]
\[
\leq J_t^\alpha \| u_k^{i+p}-u_k^i  \|_{L_2(G)}^2 +  J_t^\alpha c_{2,2}^k(t).
\]
Following straightforward computations and applying the estimate (\ref{integral}), (\ref{u_x}) can be reformulated as:
\[
\sum\limits_{k=1}^\infty\int\limits_0^t\int\limits_G\left[| \nabla(u_k^{i+p}-u_k^i) |^2 + \lambda_k | u_k^{i+p}-u_k^i |^2\right]dx d\tau \leq  \frac{3\, T}{2}\,\max_{t\in [0,T]} \sum\limits_{k=1}^\infty \| u_k^{i+p}-u_k^i \|_{L_2(G)}^2
\]
\[
+\frac{T\,C_\varepsilon\,  f_0^2}{2 } \, ||\nabla \omega||_{L_2(\Omega)}^2 \max_{x\in G} \int\limits_\Omega |f(t,x,y)|^2 dy\,\sum_{j=1}^\infty \lambda_j^{2 \tau} ||u_j^{i+p} - u_j^i||_{L_2(G)}^2, \quad t\in [0,T]. 
\]
The estimate (\ref{fundamental}) indicates that the sequences on the left-hand side are fundamental. Thus, for every \(0 < t \leq T\), as \(i \to \infty\), it follows that:
\begin{equation}\label{converUx}
\begin{cases}
\int\limits_0^t \| \nabla u_k^i(\tau, \cdot) - \nabla u_k(\tau, \cdot) \|_{L_2(G)}^2 d\tau \to 0, , \quad k=1,2, 
\cdots,\\
\int\limits_0^t \lambda_k\| u_k^i(\tau, \cdot) - u_k(\tau, \cdot) \|_{L_2(G)}^2 d\tau\to 0, , \quad k=1,2, 
\cdots,\\
\int\limits_0^t\sum_{k=1}^{\infty}  \| \nabla u_k^i (\tau, \cdot) - \nabla  u_k (\tau, \cdot) \|_{L_2(G)}^2 d\tau \to 0, \\
\int\limits_0^t\sum_{k=1}^{\infty} \lambda_k \| u_k^i(\tau, \cdot) - u_k(\tau, \cdot) \|_{L_2(G)}^2 d\tau\to 0.
\end{cases}
\end{equation}

Let us set $w = D^\alpha_t(u_k^{i+1} - u_k^i)$ in equality (\ref{Dalpha}) and replicate the reasoning provided in the proof of (\ref{Dalpha1}). Then, employing arguments analogous to those in the proof of (\ref{fundamental}), we confirm that the sequence $\sum_{k=1}^\infty \int_0^t \int_G | D_t^\alpha u_k^i |^2 \, dx \, d\tau$ is fundamental. Consequently, for every \(0 < t \leq T\), as \(i \to \infty\), it follows that:
\begin{equation}\label{fundamental2}
\begin{cases}
\int\limits_0^t \| D_t^\alpha u_k^i (\tau, \cdot) - D_t^\alpha u_k(\tau, \cdot) \|_{L_2(G)}^2 d\tau \to 0, \quad k=1,2, 
\cdots, \\
    \int\limits_0^t\sum_{k=1}^{\infty}  \| D_t^\alpha u_k^i (\tau, \cdot) - D_t^\alpha u_k(\tau, \cdot) \|_{L_2(G)}^2 d\tau \to 0.
    \end{cases}
\end{equation}

It is easy to verify that the right-hand side of equation (\ref{9}) converges in \(L_2(G)\) as $i \to \infty$. Indeed (see (\ref{C})), first note
\begin{equation}\label{right part}
\left| \sum_{k=1}^\infty (\nabla v_k, \nabla \omega) (u^{i}_k(t,x) - u_k(t,x) )\right|^2 \,  \leq C_\varepsilon ||\nabla \omega||_{L_2(\Omega)}^2 \sum_{k=1}^\infty \lambda_k^{2 \tau} |u_k^{i}(t,x) - u_k(t,x)|^2.
\end{equation}
Now, we integrate (\ref{right part}) over the domain $G$ and apply the relation (\ref{converU}). Then we obtain the required convergence.

\section{Proof of Theorem \ref{theorem 1}}
Let us pass to the limit $i\to \infty$ in (\ref{18}). Then, by virtue of the convergence (\ref{converU}) and Parseval's equality, we see that for each $t\in [0,T]$ the function
\[
u (t, x, y) = \sum\limits_{k=1}^{\infty} u_k (t, x) v_k (y),
\]
is well defined element of $L_2(D)$ and (see Lemma \ref{u_k})
\begin{equation}\label{uestimate}
\max_{[0,T]}\| u(t, \cdot, \cdot) \|^2_{\tau, G}\leq 2 A_0.
\end{equation}
Similarly, if we pass to the limit $i\to \infty$ in (\ref{u_x}), then by virtue of the convergence of (\ref{converUx}) we will obtain
 \[
\int\limits_0^t\left[||\nabla_x u(\tau. \cdot, \cdot) ||_{L_2(D)}^2 +||\nabla_y u(\tau. \cdot, \cdot) ||_{L_2(D)}^2\right]d\tau \leq \frac{\Gamma(\alpha)T^{1-\alpha}}{2}||\varphi||^2_{L_2(D)}+ 3 A_0 T +\frac{T}{2} \, A_1 
\]
\begin{equation}\label{u_xu_y}
+T\, A_0 C_\varepsilon\, f_0^2\, ||\nabla \omega||_{L_2(\Omega)}^2 \max_{x\in G} ||f(t,x,\cdot)||_{L_2(\Omega)}^2, \quad t\in [0,T]. 
\end{equation}
For any $t\in [0,T]$ the estimate (\ref{Dalpha1}) implies  
\begin{equation}\label{DalphaEstimate}
\int\limits_0^t|| D_t^\alpha u(\tau, \cdot, \cdot)||_{L_2(D)}^2d\tau \leq \frac{\Gamma(\alpha)T^{1-\alpha}}{2} \left(\| \nabla_x\varphi \|^2_{L_2(D)} +  \| \nabla_y\varphi \|^2_{L_2(D)}\right)
\end{equation}
\[
+T\max_t\left[f_0 ^2(\psi_0 + g_0)^2 \|f(t, \cdot, \cdot)\|^2_{L_2(D)}
+ \|g(t, \cdot, \cdot)\|^2_{L_2(D)}\right]
\]
\[
+2\, A_0 C_\varepsilon\, f_0^2\, ||\nabla \omega||_{L_2(\Omega)}^2 \max_{x\in G} ||f(t,x,\cdot)||_{L_2(\Omega)}^2.
\]

 \begin{лем}\label{hL2} For each $t\in [0,T]$, the function \( h (t, x) \) defined by the formula \eqref{7} is an element of $L_2(G)$ and the estimate
 \[
\int_G|h(t,x)|^2dx \leq 4 f_0^2\left(|G|\,\psi_0^2 +  2\, A_0\, C_\varepsilon ||\nabla \omega||_{L_2(\Omega)}^2\right)+ 2 |G|\,g_0^2, \,\, t\in [0,T],
\]
holds.
\end{лем}
\begin{proof} Due to inequality $(a+b +c + d)^2 \leq 4 (a^2+  b^2 + c^2 +d^2)$, it follows from definition (\ref{h_final}) of $h(t,x)$
  \begin{equation}\label{h1}
|h(t,x)|^2\leq 4 f_0^2 \left( |D_t^{\alpha} \psi(t, x)|^2 + |\Delta_x \psi(t, x)|^2 + |(g(t, x, \cdot), \omega)|^2 + \left|\sum_{j=1}^\infty (\nabla v_j, \nabla \omega) u_j(t,x) \right|^2 \right).
\end{equation}
For the sum we have the estimate (see (\ref{C}))
\[
\left| \sum_{j=1}^\infty (\nabla v_j, \nabla \omega) u_j(t,x) \right|^2 \leq  C_\varepsilon ||\nabla \omega||_{L_2(\Omega)}^2 \sum_{j=1}^\infty \lambda_j^{2 \tau} |u_j^{n-1}(t,x)|^2; , \quad x \in G, \quad t \geq 0.
\]
If we integrate this inequality over $x\in G$, then by virtue of estimate (\ref{18}) we will have
\begin{equation}\label{h2}
\int_G\left| \sum_{j=1}^\infty (\nabla v_j, \nabla \omega) u_j(t,x) \right|^2 dx \leq 2\, A_0\, C_\varepsilon ||\nabla \omega||_{L_2(\Omega)}^2.
\end{equation}
Now we integrate (\ref{h1}) with respect to $x \in G$. Then by (\ref{h2}), we get the required estimate.  
\end{proof}

Now, returning to the proof of Theorem \ref{theorem 1}, we first show the existence of a weak solution of the inverse problem. We integrate \eqref{9} over $t\in [0,T]$, and then pass to the limit $i \rightarrow \infty$ and, taking into account \eqref{converU} - \eqref{right part}, we obtain the following equalities, valid for an arbitrary \( w \in \dot{W}_2^1(0, 1) \) and for all $k$: \begin{equation}\label{22}
\int_0^t\left[( D_t^{\alpha} u_k, w ) + (\nabla u_k, \nabla w ) + \lambda_k ( u_k, w )\right] d\tau =\int_0^t\left[( g_k, w )+( f_kh, w )\right] d\tau.  
\end{equation}
Next, we take the function $w=v(x,y) v_k(y)$, $v\in \dot{W}_2^1(\Omega)$, and integrate equality (\ref{22}) over $y\in \Omega$. Then, summing over $k$ from 1 to $\infty$, we obtain for all $t\in [0, T]$ (see (\ref{gradientlambda}))
 \begin{equation}\label{solution_t}
\int_0^t \int_D \left[D_t^\alpha u v + \nabla_x u \nabla_x v +\nabla_y u \nabla_y v \right]\, dx \, dy d\tau = \int_0^t\int_D \left[f h v + g v \right]\, dx \, dy d\tau.
\end{equation}

Note that if $s(t)$ is integrable in any subset $(0,t)$ of the interval $[0,T]$ and $\int_0^t s(\tau) d\tau =0$, then obviously $s(t)=0$ almost everywhere on $[0,T]$. Therefore, the equality (\ref{solution_t}) coincides with (\ref{equation}). Consequently, the function $u(t,x,y)$ defined by formula (\ref{6}), together with $h$, defined by formula (\ref{h_final}), is a weak solution of the inverse problem, i.e., it satisfies all the conditions of definition \ref{opr1}. As for the estimates \eqref{uestimate1} - (\ref{DalphaInt}), they are proved above (see the estimates (\ref{uestimate}) - (\ref{DalphaEstimate})).  The estimate of $h$ is established in Lemma \ref{hL2}.

Let us prove the uniqueness of the weak solution of the inverse problem.

Suppose the opposite, i.e. there are two solutions to the inverse problem: \( (u_1, h_1) \) and \( (u_2, h_2) \).
Then the functions \( u = u_1 - u_2 \) and \( h = h_1 - h_2 \) satisfy all the conditions of Definition \ref{opr1} with the functions
\( \varphi(x,y) \equiv 0 \), \( g(t,x,y) \equiv 0 \). Let us denote
\[
u_k(t, x) = \int\limits_\Omega u(t, x, y) v_k (y) \, dy.
\]
Let us consider the sum that defines the function $h(t,x)$ (see (\ref{h_final})):
\[
\sum_{j=1}^\infty (\nabla v_j, \nabla \omega) u_j(t,x).
\]
Estimate (\ref{h2}) means that this function
exists almost everywhere in $G$ and for all $t\in [0, T]$. Therefore, the function
\[
h(t,x) = - \frac{\sum_{j=1}^\infty (\nabla v_j, \nabla \omega) u_j(t,x)}{(f(t, x, \cdot), \omega)}
\]
(see (\ref{h_final})) is correctly defined. Then for functions $u_k$ one has (see \eqref{12}):

\[
\frac{1}{2} D_t^\alpha \| u_k \|^2_{L_2(G)} + \|\nabla u_k \|^2_{L_2(G)} + \lambda_k \| u_k \|^2_{L_2(G)}\leq 
\]
\[
\leq \left| \int_G u_k(t,x)\frac{f_k (t,x)}{(f(t, x, \cdot), \omega)} \sum_{j=1}^\infty (\nabla v_j, \nabla \omega) u_j(t,x) \,dx \right|
\]
Since $M_k (t,x) \equiv 0$, then, repeating the same reasoning as in the proof \eqref{18}, we get:
\[
\max_{[0,T]}\sum_{k=1}^\infty \lambda_k^{2 \tau} \| u_k \|_{L_2(G)}^2  \leq 0
\]
(note that $A_0= 0$). Hence,
\[
u_k (t,x) = 0 \quad \text{a.e. in } G \text{ for all  } k \geq 1 \text{ and } t\in [0, T].
\]
 Therefore,  $h(t,x) = 0$ \text{a.e. in} $G$ for all $t\in [0, T]$. Completeness of the system $\{ v_k (y) \}$ implies
\[
u(t,x,y) = 0 \quad \text{for almost all}\,\,  x\in G, \, y \in \Omega \text{ and } t\in [0, T].
\]
Thus, the theorem is completely proved.

\begin{зам}\label{zam.3} It is important to emphasize that theorem \ref{theorem 1} remains valid for parabolic equations as well. In this case, instead of Alikhanov's estimate, one should use the usual equality $\frac{d}{dt} u^2= 2 u \frac{d}{dt} u$.
\end{зам}

\section{Conclusions}
 
In this paper, we study a new inverse problem for subdiffusion equations, namely, the problem of determining the coefficient of the source function. The elliptic part of the equation is given by $\Delta_x u + \Delta_y u$, $x\in \mathbb{R}^m$, $y\in \mathbb{R}^n$, and the unknown function depends on both time and a part of the spatial variables, denoted by $h(t, x)$.
Moreover, the overdetermination condition is given in integral form. To our knowledge, such a problem has not been studied before for subdiffusion equations or even for diffusion equations. The existence and uniqueness of a weak solution to the inverse problem are established, as well as coercive estimates.

The choice of constructing a weak solution is motivated by the relative simplicity of obtaining a priori estimates in such a formulation. If we consider the classical formulation of the problem (\ref{9}), (\ref{10}) (for example, see \cite{AMux2}), then we can prove the existence of a classical solution to the inverse problem. However, this will be the subject of a future paper.

\

\begin{center}
ACKNOWLEDGEMENTS    
\end{center}

 The authors acknowledge financial support from the Ministry of Innovative Development of the Republic of Uzbekistan, Grant No F-FA-2021-424.


\end{document}